\documentclass[12pt]{article}
\usepackage{makeidx}
\usepackage{graphicx}
\usepackage{amsmath}
\usepackage{amsfonts}
\usepackage{amssymb}
\newtheorem{theorem}{Theorem}

\newtheorem{definition}[theorem]{Definition}
\newtheorem{example}[theorem]{Example}

\newtheorem{proposition}[theorem]{Proposition}

\newenvironment{proof}[1][Proof]{\textbf{#1.} }{\ \rule{0.5em}{0.5em}}

\begin{document}

\title{Pad\'{e} approximation for a multivariate Markov transform }
\author{Ognyan Kounchev and Hermann Render}
\maketitle
\begin{abstract}
Methods of Pad\'{e} approximation are used to analyse a multivariate Markov
transform which has been recently introduced by the authors. The first main
result is a characterization of the rationality of the Markov transform via
Hankel determinants. The second main result is a cubature formula for a
special class of measures.

\textrm{Acknowledgement}: \emph{The authors thank the Alexander von Humboldt
Foundation for support in the framework of the Institutes partnership project
and the Feodor-Lynen programme. The second author is supported in part by
Grant MTM2006-13000-C03-03 of the D.G.I. of Spain.}


\end{abstract}

\section{Introduction}

Let $\sigma$ be a non-negative finite measure on a subinterval $\left[
a,b\right]  $ of the real line $\mathbb{R.}$ Then the numbers $\int_{a}%
^{b}x^{l}d\sigma\left(  x\right)  $ are called the \emph{moments of the
measure }$\sigma.$ The \emph{Markov transform} of $\sigma$ is defined for
$\zeta\in\mathbb{C}\setminus\left[  a,b\right]  $ by the formula
\begin{equation}
\widehat{\sigma}\left(  \zeta\right)  :=\int_{a}^{b}\frac{1}{\zeta-x}%
d\sigma\left(  x\right)  \text{.}\label{defMarkov}%
\end{equation}
In the theory of moments Pad\'{e} approximation of the Markov transform
$\widehat{\sigma}\left(  \zeta\right)  $ is an important tool, see
\cite{Akhi65}, \cite{Brez80}, \cite{Brez91} or \cite{NiSo91} and section 6.
Here Pad\'{e} approximation is performed at the point $\infty$, so we consider
the asymptotic expansion
\begin{equation}
\widehat{\sigma}\left(  \zeta\right)  =\sum_{l=0}^{\infty}\int_{a}^{b}%
x^{l}d\sigma\left(  x\right)  \frac{1}{\zeta^{l+1}}\text{ for }\left|
\zeta\right|  >R.\label{asympst}%
\end{equation}

Let now $\mu$ be a signed measure $\mu$ on the euclidean space $\mathbb{R}%
^{d}$ with support in the closed ball $\overline{B_{R}}:=\left\{
x\in\mathbb{R}^{d}:\left|  x\right|  \leq R\right\}  $ where $\left|
x\right|  :=\sqrt{x_{1}^{2}+...+x_{d}^{2}}$ is the euclidean distance for
$x=\left(  x_{1},...,x_{d}\right)  \in\mathbb{R}^{d}$. In \cite{KoReHirosh} we
introduced a \emph{multivariate Markov transform} for the measure $\mu$ by the
formula
\begin{equation}
\widehat{\mu}\left(  \zeta,\theta\right)  =\int_{\mathbb{R}^{d}}\frac
{\zeta^{d-1}}{r\left(  \zeta\theta-x\right)  ^{d}}d\mu\left(  x\right)  \text{
for }\left|  \zeta\right|  >R,\theta\in\mathbb{S}^{d-1}.\label{defmst}%
\end{equation}
Here $\mathbb{S}^{d-1}:=\left\{  x\in\mathbb{R}^{d}:\left|  x\right|
=1\right\}  $ is the unit sphere, and $\zeta$ is a complex number with
$\left|  \zeta\right|  >R.$ In the denominator, the expression $r\left(
\zeta\theta-x\right)  $ is the analytic continuation of the function
$\rho\longmapsto\left|  \rho\theta-x\right|  $ defined for $\rho\in\mathbb{R}$
with $\rho>R$, see Section 3 for details. The motivation for this definition
stems from the work of \emph{N. Aronszajn} about polyharmonic functions and
the work of \emph{L.K. Hua} about harmonic analysis on Lie groups, see
\cite{ACL83}, \cite{Hua63} or \cite{KoReHirosh}. Following the analogy with
the one-dimensional case, we consider the asymptotic expansion of the
multivariate Markov transform. From the growth behaviour at infinity of the
kernel $\zeta^{d-1}/r\left(  \zeta\theta-x\right)  ^{d}$ it is easily seen
that the asymptotic expansion is of the form
\begin{equation}
\widehat{\mu}\left(  \zeta,\theta\right)  =\sum_{l=0}^{\infty}f_{l}\left(
\theta\right)  \frac{1}{\zeta^{l+1}}\label{eqmuasym}%
\end{equation}
for $\left|  \zeta\right|  >R$ and $\theta\in\mathbb{S}^{d-1}$ where $f_{l}: $
$\mathbb{S}^{d-1}\rightarrow\mathbb{C}$ are continuous functions. The aim of
this paper is to show that methods from Pad\'{e} approximation can be
successfully used for an analysis of the multivariate Markov transform.
Roughly speaking, we shall perform in (\ref{eqmuasym}) the classical
univariate Pad\'{e} approximation for each fixed $\theta\in\mathbb{S}^{d-1} $
obtaining a Pad\'{e} pair $(Q_{n}\left(  \zeta,\theta\right)  ,P_{n}\left(
\zeta,\theta\right)  )$.

Let us describe the results in the paper: In section 2 we shall first review
the basic notions from Pad\'{e} approximation which are needed in the paper.
In section 3 the asymptotic expansion defined in (\ref{eqmuasym}) will be
investigated. It turns out that each coefficient function $f_{l}$ in
(\ref{eqmuasym}) is a finite sum of spherical harmonics of degree $\leq l,$
and each $f_{l}$ is the restriction of a homogeneous polynomial $F_{l}\left(
x\right)  $ of degree $l$ to the unit sphere. The \emph{Hankel determinant} of
the multivariate Markov transform $\widehat{\mu}$ (or a measure $\mu)$ is
defined by the expression
\begin{equation}
H_{n}\left(  \mu,\theta\right)  :=\det\left(
\begin{array}
[c]{llll}%
f_{0}\left(  \theta\right)  & f_{1}\left(  \theta\right)  & \cdots &
f_{n-1}\left(  \theta\right) \\
f_{1}\left(  \theta\right)  & f_{2}\left(  \theta\right)  & \cdots &
f_{n}\left(  \theta\right) \\
\cdots & \cdots & \cdots & \cdots\\
f_{n-1}\left(  \theta\right)  & f_{n}\left(  \theta\right)  & \cdots &
f_{2n-2}\left(  \theta\right)
\end{array}
\right)  .\label{defhn}%
\end{equation}
In section 4 we show that the Hankel determinant $H_{n}\left(  \mu
,\theta\right)  $ of a measure $\mu$ is the restriction of a homogeneous
polynomial of degree $n\left(  n-1\right)  $ to the unit sphere. In section 5
we shall prove a Kronecker type theorem: the Hankel determinants $H_{n}\left(
\mu,\theta\right)  $ are zero for all large $n$ if and only if the function
$\zeta\longmapsto\widehat{\mu}\left(  \zeta,\theta\right)  $ is rational for
each $\theta\in\mathbb{S}^{d-1}.$ Moreover, this is equivalent to the
rationality of the multivariate Markov transform $\widehat{\mu}.$

A measure $\mu$ is called \emph{Hankel positive} if the Hankel determinants
$H_{n}\left(  \mu,\theta\right)  $ are strictly positive for all natural
numbers $n$ and for all $\theta\in\mathbb{S}^{d-1}.$ In section 6 we prove
that for each Hankel positive measure $\mu$ there exists a non-negative
measure $\mu_{n}$ which is equal to $\mu$ for all polynomials of degree
$\leq2n-1$ and which has support contained in an algebraic variety. Further we
characterize Hankel positivity by an extension property of the multivariate
Markov transform.

Finally we need some notations from harmonic analysis. A function
$Y:\mathbb{S}^{d-1}\rightarrow\mathbb{C}$ is called a \emph{spherical
harmonic} of degree $k\in\mathbb{N}_{0}$ if there exists a \emph{homogeneous}
\emph{harmonic} polynomial $P\left(  x\right)  $ of degree $k$ (in general,
with complex coefficients) such that $P\left(  \theta\right)  =Y\left(
\theta\right)  $ for all $\theta\in\mathbb{S}^{d-1}.$ Throughout the paper we
assume that $Y_{k,m}\left(  x\right)  $, $m=1,...,a_{k},$\ is a basis of the
set of all harmonic homogeneous polynomials of degree $k$\ which are
orthonormal with respect to scalar product
\[
\left\langle f,g\right\rangle _{\mathbb{S}^{d-1}}:=\int_{\mathbb{S}^{d-1}%
}f\left(  \theta\right)  \overline{g\left(  \theta\right)  }d\theta.
\]
Here $a_{k}$ denotes the dimension of the space of all harmonic homogeneous
polynomials of degree $k.$ By $\omega_{d}$ we denote the surface area of
$\mathbb{S}^{d-1}.$

\section{Basic facts from Pad\'{e} Approximation}

At first let us recall some basic facts from Pad\'{e} Approximation (we refer
to \cite{NiSo91} for proofs): let $f$ be a holomorphic function for $\zeta
\in\mathbb{C},$ $\left|  \zeta\right|  >R,$ of the form
\[
f\left(  \zeta\right)  =\sum_{l=0}^{\infty}f_{l}\frac{1}{\zeta^{l+1}}.
\]
Let $n$ be a natural number. Then there exists a polynomial $P_{n}\neq0$ of
degree $\leq n$ such that
\begin{equation}
P_{n}\left(  \zeta\right)  f\left(  \zeta\right)  -Q_{n}\left(  \zeta\right)
=\sum_{l=n}^{\infty}f_{l}\frac{1}{\zeta^{l+1}}\label{eqPade}%
\end{equation}
where $Q_{n}$ is the polynomial part of the series $P_{n}\left(  z\right)
f\left(  \zeta\right)  ;$ it is easy to see that $Q_{n}$ has degree $\leq
n-1$. A pair $(P_{n},Q_{n})$ is called an $n$\emph{-th Pad\'{e} pair} if
$P_{n}$ and $Q_{n}$ are polynomials, $P_{n}\neq0,$ $\deg P_{n}\leq n$ and
$\deg Q_{n}\leq n-1,$ and they satisfy (\ref{eqPade}). An index $n$ is called
\emph{normal }if for any $n$-th Pad\'{e} pair $(P_{n},Q_{n}$) the polynomial
$\zeta\mapsto P_{n}\left(  \zeta\right)  $ has degree exactly $n.$ Proposition
3.2 in \cite{NiSo91} shows that $n$ is normal if and only if the \emph{Hankel
determinant}
\begin{equation}
H_{n}\left(  f\right)  :=\det\left(
\begin{array}
[c]{llll}%
f_{0} & f_{1} & \cdots &  f_{n-1}\\
f_{1} & f_{2} & \cdots &  f_{n}\\
\cdots & \cdots & \cdots & \cdots\\
f_{n-1} & f_{n} & \cdots &  f_{2n-2}%
\end{array}
\right) \label{defhn2}%
\end{equation}
is not zero. If $n$ is normal then the polynomial
\[
P_{n}\left(  \zeta\right)  :=\det\left(
\begin{array}
[c]{llll}%
f_{0} & f_{1} & \cdots &  f_{n}\\
\cdots & \cdots & \cdots & \cdots\\
f_{n-1} & \cdots & \cdots &  f_{2n-1}\\
1 & \zeta & \cdots & \zeta^{n}%
\end{array}
\right)
\]
has exact degree $n$ and $(P_{n},Q_{n})$ is an $n$-th Pad\'{e} pair where
$Q_{n}$ is the polynomial part of $P_{n}\left(  z\right)  f\left(  z\right)
.$ For arbitrary $n,$ the rational function
\[
\pi_{n}\left(  \zeta\right)  :=\frac{Q_{n}\left(  \zeta\right)  }{P_{n}\left(
\zeta\right)  }%
\]
is called the $n$\emph{-th diagonal Pad\'{e} approximant} of $f$.

\section{Asymptotic expansion of the multivariate Markov transform}

Following the exposition in \cite[Section 2.2]{ACL83} we show that the
multivariate Markov transform is well-defined. Let us set $r\left(  x\right)
:=\left|  x\right|  .$ For $\rho>0$ and $\theta\in\mathbb{S}^{d-1}$ and
$x=\left(  x_{1},...,x_{d}\right)  $ we have $r^{2}\left(  \rho\theta
-x\right)  =\rho^{2}-2\rho\left\langle \theta,x\right\rangle +\left|
x\right|  ^{2}$ where $\left\langle \theta,x\right\rangle $ is the usual inner
product in $\mathbb{R}^{d}.$ We replace $\rho$ by a complex number $\zeta$ and
obtain
\[
r^{2}\left(  \zeta\theta-x\right)  =\zeta^{2}-2\zeta\left\langle
\theta,x\right\rangle +\left|  x\right|  ^{2}=\left(  \zeta-\left\langle
\theta,x\right\rangle \right)  ^{2}+\left|  x\right|  ^{2}-\left|
\left\langle \theta,x\right\rangle \right|  ^{2}.
\]
Note that $\left|  x\right|  ^{2}-\left|  \left\langle \theta,x\right\rangle
\right|  ^{2}\geq0$ for each $\theta\in\mathbb{S}^{d-1}.$ If we define
\[
a\left(  \theta,x\right)  :=\left\langle \theta,x\right\rangle +i\sqrt{\left|
x\right|  ^{2}-\left|  \left\langle \theta,x\right\rangle \right|  ^{2}}%
\]
then $r^{2}\left(  \zeta\theta-x\right)  =\left(  \zeta-a\left(
\theta,x\right)  \right)  \left(  \zeta-\overline{a\left(  \theta,x\right)
}\right)  . $ Since $\left|  a\left(  \theta,x\right)  \right|  ^{2}=\left|
x\right|  ^{2}$ it follows that $r^{2}\left(  \zeta\theta-x\right)  \neq0$ for
all $\left|  \zeta\right|  >\left|  x\right|  .$ Next we see that the function
$g,$ defined by
\[
g\left(  \zeta\right)  :=\frac{r^{2}\left(  \zeta\theta-x\right)  }{\zeta^{2}%
}=\left(  1-\frac{a\left(  \theta,x\right)  }{\zeta}\right)  \left(
1-\frac{\overline{a\left(  \theta,x\right)  }}{\zeta}\right)
\]
for $\left|  \zeta\right|  >\left|  x\right|  ,$ has the property that
$g\left(  \zeta\right)  \notin\left(  -\infty,0\right]  $: since $\left|
\frac{a\left(  \theta,x\right)  }{\zeta}\right|  <1$ and $\left|
\frac{\overline{a\left(  \theta,x\right)  }}{\zeta}\right|  <1$ it follows
that $1-\frac{a\left(  \theta,x\right)  }{\zeta}$ and $1-\frac{\overline
{a\left(  \theta,x\right)  }}{\zeta}$ are in the right half plane, i.e. that
their real parts are strictly positive, and therefore $g\left(  \zeta\right)
\notin\left(  -\infty,0\right]  .$ Using the square root function $\sqrt
{\cdot}$ defined on $\mathbb{C}\setminus\left(  -\infty,0\right]  $ one can
define for $\left|  \zeta\right|  >\left|  x\right|  $ the analytic function
\[
\zeta\longmapsto\sqrt{\frac{r^{2}\left(  \zeta\theta-x\right)  }{\zeta^{2}}}.
\]
It follows from these facts that the multivariate Markov transform
$\widehat{\mu}\left(  \zeta,\theta\right)  $ is well-defined.

Further we will make use of a real version of the multivariate Markov
transform which we define by (note that we use $d$ instead of $d-1$ as
exponent in the nominator)
\begin{equation}
\widehat{\mu}_{\text{real}}\left(  y\right)  :=\int_{\mathbb{R}^{d}}%
\frac{\left|  y\right|  ^{d}}{r\left(  y-x\right)  ^{d}}d\mu\left(  x\right)
\text{ for }y\in\mathbb{R}^{d}\text{ with }\left|  y\right|
>R.\label{muschlange}%
\end{equation}
The \emph{real Markov transform} $\widehat{\mu}_{\text{real}}\left(  y\right)
$ is related to $\widehat{\mu}\left(  \zeta,\theta\right)  $ in the following
way: using results about harmonicity hulls and Lie norms (see \cite[p.
64]{ACL83}) one may show that the function $y\longmapsto\widehat{\mu
}_{\text{real}}\left(  y\right)  $ has an holomorphic extension to a natural
set $C_{R}$ in the complex space $\mathbb{C}^{d},$ and the extension will be
denoted by $\widehat{\mu}_{\text{real}}\left(  z\right)  $ for complex $z\in
C_{R}.$ The set $C_{R}$ is the set of all $z=\left(  z_{1},....,z_{d}\right)
\in\mathbb{C}^{d}$ such that
\begin{equation}
L_{-}\left(  z\right)  :=\sqrt{\left|  z\right|  ^{2}-\sqrt{\left|  z\right|
^{4}-\left|  q\left(  z\right)  \right|  ^{2}}}>R\label{defLminus}%
\end{equation}
where we have defined $\left|  z\right|  ^{2}=\left|  z_{1}\right|
^{2}+...+\left|  z_{d}\right|  ^{2}$ and $q\left(  z\right)  =z_{1}%
^{2}+...+z_{d}^{2}.$ The set $C_{R}$ is connected and open, and it contains
all points $\zeta\cdot\theta$ with $\zeta\in\mathbb{C}$, $\left|
\zeta\right|  >R$ and $\theta\in\mathbb{S}^{d-1}$. The Markov transforms
$\widehat{\mu}\left(  \zeta,\theta\right)  $ and $\widehat{\mu}_{\text{real}%
}\left(  z\right)  $ are related by the simple formula
\begin{equation}
\widehat{\mu}_{\text{real}}\left(  \zeta\theta\right)  =\zeta\widehat{\mu
}\left(  \zeta,\theta\right)  \text{ for all }\zeta\in\mathbb{C},\left|
\zeta\right|  >R\text{ and }\theta\in\mathbb{S}^{d-1}.\label{eqmureell}%
\end{equation}

Next we want to describe the asymptotic expansion of the multivariate Markov
transform. Using the \emph{Gau\ss\ decomposition} of a polynomial (see Theorem
5.5 in \cite{ABR92}, \cite{StWe71}, or \cite{Koun00}) it is easy to see that
the system
\[
\left|  x\right|  ^{2s}Y_{k,m}\left(  x\right)  ,s,k\in\mathbb{N}%
_{0},m=1,...,a_{k}%
\]
is a basis of the set of all polynomials. The numbers
\begin{equation}
c_{s,k,m}:=\int_{\mathbb{R}^{d}}\left|  x\right|  ^{2s}\overline
{Y_{k,m}\left(  x\right)  }d\mu\left(  x\right)  ,\quad s,k\in\mathbb{N}%
_{0},m=1,...,a_{k}\label{distributed}%
\end{equation}
are sometimes called the \emph{distributed moments}, see \cite{kounchev87}.
For a treatment and formulation of the \emph{multivariate moment problem} we
refer to \cite{Fugl83} and \cite{ACL83}. From \cite{KoReHirosh} we cite

\begin{theorem}
\label{T3}Let $\mu$ be a signed measure on $\mathbb{R}^{d}$ with support in
the closed ball $\overline{B_{R}}$. Then for all $\left|  \zeta\right|  >R$
and for all $\theta\in\mathbb{S}^{d-1}$ the following relation holds
\begin{equation}
\widehat{\mu}\left(  \zeta,\theta\right)  =\sum_{s=0}^{\infty}\sum
_{k=0}^{\infty}\sum_{m=1}^{a_{k}}\frac{Y_{k,m}\left(  \theta\right)  }%
{\zeta^{2s+k+1}}\int_{\mathbb{R}^{d}}\left|  x\right|  ^{2s}\overline
{Y_{k,m}\left(  x\right)  }d\mu\left(  x\right)  .\label{muhatrep2}%
\end{equation}
\end{theorem}

For $f_{l}$ defined in (\ref{eqmuasym}), a rearrangement of the series
(\ref{muhatrep2}) in powers $\zeta^{l+1}$ yields the relation
\begin{equation}
f_{l}\left(  \theta\right)  =\sum_{t=0}^{\left[  \frac{l}{2}\right]  }%
\sum_{m=1}^{\alpha_{l-2t}}c_{t,l-2t,m}Y_{l-2t,m}\left(  \theta\right)
,\label{deffl}%
\end{equation}
where $\left[  x\right]  $ denotes the largest integer $n$ such that $n\leq x.$

\begin{proposition}
\label{P8}For each $l\in\mathbb{N}_{0}$ the coefficient function $f_{l}$ in
(\ref{eqmuasym}) is a finite sum of spherical harmonics of degree $\leq l$.
Moreover, there exists a homogeneous polynomial $F_{l}\left(  x\right)  $ of
degree $l$ such that
\[
F_{l}\left(  \zeta\theta\right)  =\zeta^{l}f_{l}\left(  \theta\right)  \text{
for all }\theta\in\mathbb{S}^{d-1}\text{ and }\zeta\in\mathbb{C}.
\]
\end{proposition}%

\proof
Formula (\ref{deffl}) shows that $f_{l}$ is a sum of sphercial harmonics of
degree $\leq l.$ Define a homogeneous polynomial $F_{l}$ of degree $l$ by
\begin{equation}
F_{l}\left(  x\right)  :=\sum_{t=0}^{\left[  \frac{l}{2}\right]  }\sum
_{m=1}^{\alpha_{l-2t}}c_{t,l-2t,m}\left|  x\right|  ^{2t}Y_{l-2t,m}\left(
x\right)  .\label{defpl}%
\end{equation}
By inserting $x=\rho\theta$ in (\ref{defpl}) for positive $\rho$ we obtain
$F_{l}\left(  \rho\theta\right)  =\rho^{l}f_{l}\left(  \theta\right)  .$ Since
$\rho\longmapsto F_{l}\left(  \rho\theta\right)  $ is holomorphic we may
replace $\rho$ by a complex number $\zeta.$ The proof is finished.%
\endproof

The coefficient function $f_{l}$ can also be described by Legendre polynomials
$P_{k}\left(  t\right)  $ of degree $k$ and dimension $d,$ for definition see
\cite{Mull66}. Clearly (\ref{deffl}) and (\ref{distributed}) implies that
\[
f_{l}\left(  \theta\right)  =\int_{\mathbb{R}^{d}}\sum_{t=0}^{\left[
l/2\right]  }\left|  x\right|  ^{2t}\sum_{m=1}^{a_{l-2t}}Y_{l-2t,m}\left(
x\right)  \cdot Y_{l-2t,m}\left(  \theta\right)  d\mu\left(  x\right)  .
\]
The addition theorem for spherical harmonics (see \cite{Mull66}) says that
\[
\sum_{m=1}^{a_{k}}Y_{k,m}\left(  x\right)  \cdot Y_{k,m}\left(  \theta\right)
=\left|  x\right|  ^{k}a_{k}P_{k}\left(  \left\langle \frac{x}{\left|
x\right|  },\theta\right\rangle \right)  ,
\]
so one obtains the alternative description
\[
f_{l}\left(  \theta\right)  =\sum_{t=0}^{\left[  l/2\right]  }a_{l-2t}%
\int_{\mathbb{R}^{d}}\left|  x\right|  ^{l}P_{l-2t}\left(  \left\langle
\frac{x}{\left|  x\right|  },\theta\right\rangle \right)  d\mu\left(
x\right)  .
\]
We conclude this section with some examples and results illustrating the definitions.

\begin{example}
\label{Ex0}Let $\sigma$ be a finite non-negative measure on an interval
$\left[  a,b\right]  $ with $a\geq0$ and consider the measure $\mu
=\sigma\otimes d\theta,$ i.e. for every continuous function $f$ holds
\[
\int f\left(  x\right)  d\mu:=\int_{a}^{b}\int_{\mathbb{S}^{d-1}}f\left(
r\theta\right)  d\sigma\left(  r\right)  d\theta.
\]
Then the distributed moments $c_{s,k,m}$ are zero for all $k>0$ since
$Y_{k,m}\left(  \theta\right)  $ is orthogonal to the constant function with
respect to the measure $d\theta.$ Hence (\ref{muhatrep2}) shows that
\[
\widehat{\mu}\left(  \zeta,\theta\right)  =\sum_{s=0}^{\infty}\frac{1}%
{\zeta^{2s+1}}\int_{a}^{b}r^{2s}d\sigma\left(  r\right)  =\int_{a}^{b}%
\frac{\zeta}{\zeta^{2}-r^{2}}d\sigma\left(  r\right)  .
\]
From this we conclude that for all $l\in\mathbb{N}_{0}$ and $\theta
\in\mathbb{S}^{d-1}$
\[
f_{2l}\left(  \theta\right)  =\int_{a}^{b}r^{2l}d\sigma\left(  r\right)
\text{ and }f_{2l+1}\left(  \theta\right)  =0\text{.}%
\]
\end{example}

A measure $\mu$ on {$\mathbb{R}$}$^{d}$ is called \emph{rotation invariant} if
$\mu\left(  T^{-1}\left(  B\right)  \right)  =\mu\left(  B\right)  $ for all
Borel sets $B$ and for all orthogonal linear maps $T:$ {$\mathbb{R}$}%
$^{d}\rightarrow${$\mathbb{R}$}$^{d}.$ The following result shows that a
rotation invariant measure has a Markov transform $\widehat{\mu}\left(
\zeta,\theta\right)  $ which does not depend on $\theta\in\mathbb{S}^{d-1}.$
Since the result is not needed later we omit the proof.

\begin{theorem}
\label{Rotthm}Let $\mu$ be a measure on {$\mathbb{R}$}$^{d}$ with support in
$\overline{B_{R}}$. Then $\widehat{\mu}\left(  \zeta,\theta\right)  $ is
independent of $\theta$ if and only if $\mu$ is rotation invariant. In that
case the multivariate Markov transform possesses an analytic continuation to
the upper half plane, namely
\[
\widehat{\mu}\left(  \zeta,\theta\right)  =\int\frac{\zeta}{\zeta^{2}-\left|
x\right|  ^{2}}d\mu\left(  x\right)  =\sum_{l=0}^{\infty}\int\left|  x\right|
^{2l}d\mu\frac{1}{\zeta^{2l+1}}%
\]
for all $\text{Im}\zeta>0$ and $\theta\in\mathbb{S}^{d-1}$.
\end{theorem}

\section{Multivariate Pad\'{e} approximation and Hankel determinants}

We start with the following observation:

\begin{proposition}
\label{PropHD}Let $H_{n}\left(  \mu,\theta\right)  $ be the Hankel determinant
defined in (\ref{defhn}). Then there exists a homogeneous polynomial
$\widetilde{H}_{n}\left(  x\right)  $ of degree $n\left(  n-1\right)  $ such
\[
\widetilde{H}_{n}\left(  \zeta\theta\right)  =\zeta^{n\left(  n-1\right)
}H_{n}\left(  \mu,\theta\right)  \text{ for all }\theta\in\mathbb{S}^{d-1}.
\]
\end{proposition}%

\proof
By Proposition \ref{P8} there exists a homogeneous polynomial $F_{l}\left(
x\right)  $ of degree $l$ such that $F_{l}\left(  \zeta\theta\right)
=\zeta^{l}f_{l}\left(  \theta\right)  .$ Let us define
\[
\widetilde{H}_{n}\left(  x\right)  :=\det\left(
\begin{array}
[c]{llll}%
F_{0}\left(  x\right)  & F_{1}\left(  x\right)  & \cdots &  F_{n-1}\left(
x\right) \\
F_{1}\left(  x\right)  & F_{2}\left(  x\right)  & \cdots &  F_{n}\left(
x\right) \\
\cdots & \cdots & \cdots & \cdots\\
F_{n-1}\left(  x\right)  & F_{n}\left(  x\right)  & \cdots &  F_{2n-2}\left(
x\right)
\end{array}
\right)  .
\]
Now we replace $x$ by $\zeta\theta$ and we apply the Leibniz formula for
determinants to the matrix $A=\left(  a_{i,j}\right)  _{i,j=1,...,n}$ defined
by $a_{i,j}=F_{i+j}\left(  \zeta\theta\right)  $ for $i,j=0,...,n-1$. Then
\[
\widetilde{H}_{n}\left(  \zeta\theta\right)  =\sum_{\sigma\text{ permutation
}}sign\left(  \sigma\right)  F_{0+\sigma\left(  0\right)  }\left(  \zeta
\theta\right)  ...F_{n-1+\sigma\left(  n-1\right)  }\left(  \zeta
\theta\right)  .
\]
Note that
\[
0+\sigma\left(  0\right)  +1+\sigma\left(  1\right)  +...+\left(  n-1\right)
+\sigma\left(  n-1\right)  =n\left(  n-1\right)  .
\]
It is obvious that $\widetilde{H}_{n}\left(  x\right)  $ is a homogeneous
polynomial of degree $n\left(  n-1\right)  .$ We can factor out $\zeta
^{n\left(  n-1\right)  }$ and we see that $\widetilde{H}_{n}\left(
\zeta\theta\right)  =\zeta^{n\left(  n-1\right)  } H_{n}\left(  \mu
,\theta\right)  .$
\endproof

In the following it is convenient to introduce the following notation: for a
natural number $n$ define a polynomial $\widetilde{P}_{n}\left(  \zeta
,\theta\right)  $ of a univariate variable $\zeta$ of degree $\leq n$ by
\begin{equation}
\widetilde{P}_{n}\left(  \zeta,\theta\right)  :=\det\left(
\begin{array}
[c]{llll}%
f_{0}\left(  \theta\right)  & f_{1}\left(  \theta\right)  & \cdots &
f_{n}\left(  \theta\right) \\
\cdots & \cdots & \cdots & \cdots\\
f_{n-1}\left(  \theta\right)  & \cdots & \cdots &  f_{2n-1}\left(
\theta\right) \\
1 & \zeta & \cdots & \zeta^{n}%
\end{array}
\right)  .\label{defpn}%
\end{equation}
We shall also write
\begin{equation}
\widetilde{P}_{n}\left(  \zeta,\theta\right)  =p_{0}\left(  \theta\right)
+p_{1}\left(  \theta\right)  \zeta+...+p_{n}\left(  \theta\right)  \zeta
^{n}.\label{ppointwise}%
\end{equation}
We define $\widetilde{Q}_{n}\left(  \zeta,\theta\right)  $ as the polynomial
part of $\widetilde{P}_{n}\left(  \zeta,\theta\right)  \widehat{\mu}\left(
\zeta,\theta\right)  ,$ so
\begin{equation}
\widetilde{Q}_{n}\left(  \zeta,\theta\right)  =p_{n}f_{0}\zeta^{n-1}+\left(
p_{n-1}f_{0}+p_{n}f_{1}\right)  \zeta^{n-2}+...+\left(  p_{1}f_{0}+p_{2}%
f_{1}+...+p_{n}f_{n-1}\right) \label{ppart}%
\end{equation}

From the results in section 2 the following is clear:

\begin{theorem}
\label{ThmPad}If $H_{n}\left(  \mu,\theta\right)  \neq0$\ then $(\widetilde
{P}_{n}\left(  \zeta,\theta\right)  ,\widetilde{Q}_{n}\left(  \zeta
,\theta\right)  )$\ is an $n$-th Pad\'{e} pair for the function
\[
\zeta\longmapsto\widehat{\mu}\left(  \zeta,\theta\right)  =\sum_{l=0}^{\infty
}f_{l}\left(  \theta\right)  \frac{1}{\zeta^{l+1}}%
\]
for $\left|  \zeta\right|  >R$ where $\theta\in\mathbb{S}^{d-1}$ acts as a parameter.
\end{theorem}

In Example \ref{Ex1} below we shall show that $\zeta\longmapsto\widetilde
{P}_{n}\left(  \zeta,\theta\right)  $ may be the zero polynomial for certain
$\theta\in\mathbb{S}^{d-1},$ so $(\widetilde{P}_{n}\left(  \zeta
,\theta\right)  ,\widetilde{Q}_{n}\left(  \zeta,\theta\right)  )$\emph{\ }is
\emph{not always} an $n$-th Pad\'{e} pair.

The advantage of working with $\widetilde{P}_{n}\left(  \zeta,\theta\right)  $
is seen from the following result:

\begin{theorem}
\label{ThmPoly}Let $\widetilde{P}_{n}\left(  \zeta,\theta\right)  $ and
$\widetilde{Q}_{n}\left(  \zeta,\theta\right)  $ be defined in (\ref{defpn})
and (\ref{ppart}). Then there exists a polynomial $A_{n}$ of degree $\leq
n^{2}+n$ and a polynomial $B_{n}$ of degree $\leq n^{2}+n-2$ such that
\begin{equation}
\zeta^{n^{2}}\widetilde{P}_{n}\left(  \zeta,\theta\right)  =A_{n}\left(
\zeta\theta\right)  \text{ and }\zeta^{n^{2}}\widetilde{Q}_{n}\left(
\zeta,\theta\right)  =\zeta B_{n}\left(  \zeta\theta\right)  .\label{eqPolAB}%
\end{equation}
\end{theorem}%

\proof
By Proposition \ref{P8} there exists for each $l\in\mathbb{N}_{0}$ a
homogeneous polynomial $F_{l}$ of degree $l$ such that $\zeta^{l}f_{l}\left(
\theta\right)  =F_{l}\left(  \zeta\theta\right)  .$ Let us multiply each
$j$-th column in (\ref{defpn}), $j=0,...,n,$ with $\zeta^{n-1+j}.$ Let us
define $d_{n}$ to be the sum of $n-1+j$ for $j=0,...,n.$ It follows that
$\zeta^{d_{n}}P_{n}\left(  \zeta,\theta\right)  $ is equal to
\[
\det\left(
\begin{array}
[c]{llll}%
\zeta^{n-1}f_{0}\left(  \theta\right)  & \zeta^{n}f_{1}\left(  \theta\right)
& \cdots & \zeta^{2n-1}f_{n}\left(  \theta\right) \\
\cdots & \cdots & \cdots & \cdots\\
\zeta^{n-1}f_{n-1}\left(  \theta\right)  & \cdots & \cdots & \zeta
^{2n-1}f_{2n-1}\left(  \theta\right) \\
\zeta^{n-1}1 & \zeta^{n}\zeta & \cdots & \zeta^{2n-1}\zeta^{n}%
\end{array}
\right)  .
\]
Since $F_{l}\left(  \zeta\theta\right)  =\zeta^{l}f_{l}\left(  \theta\right)
$ we obtain that $\zeta^{d_{n}}P_{n}\left(  \zeta,\theta\right)  $ is equal
to
\[
\det\left(
\begin{array}
[c]{llll}%
\zeta^{n-1}F_{0}\left(  \zeta\theta\right)  & \zeta^{n-1}F_{1}\left(
\zeta\theta\right)  & \cdots & \zeta^{n-1}F_{n}\left(  \zeta\theta\right) \\
\cdots & \cdots & \cdots & \cdots\\
F_{n-1}\left(  \zeta\theta\right)  & \cdots & \cdots &  F_{2n-1}\left(
\zeta\theta\right) \\
\zeta^{n-1}1 & \zeta^{n-1}\zeta^{2} & \cdots & \zeta^{n-1}\zeta^{2n}%
\end{array}
\right)  .
\]
From the $j$-th row factor out $\zeta^{n-1-j}$ for $j=0,...,n-1$ and from the
last one $\zeta^{n-1}.$ Then $d_{n}-\left(  n-1\right)  -\sum_{j=0}^{n-1}j $
is equal to $n^{2}.$ Hence we have proved that
\[
\zeta^{n^{2}}\widetilde{P}_{n}\left(  \zeta,\theta\right)  =\det\left(
\begin{array}
[c]{llll}%
F_{0}\left(  \zeta\theta\right)  & F_{1}\left(  \zeta\theta\right)  & \cdots
&  F_{n}\left(  \zeta\theta\right) \\
\cdots & \cdots & \cdots & \cdots\\
F_{n-1}\left(  \zeta\theta\right)  & \cdots & \cdots &  F_{2n-1}\left(
\zeta\theta\right) \\
1 & \zeta^{2} & \cdots & \zeta^{2n}%
\end{array}
\right)  .
\]
It follows that $\zeta^{n^{2}}P_{n}\left(  \zeta,\theta\right)  =A_{n}\left(
\zeta\theta\right)  $ where $A_{n}\left(  x\right)  $ is defined as
\begin{equation}
A_{n}\left(  x\right)  :=\det\left(
\begin{array}
[c]{llll}%
F_{0}\left(  x\right)  & F_{1}\left(  x\right)  & \cdots &  F_{n}\left(
x\right) \\
\cdots & \cdots & \cdots & \cdots\\
F_{n-1}\left(  x\right)  & \cdots & \cdots &  F_{2n-1}\left(  x\right) \\
1 & \left|  x\right|  ^{2} & \cdots & \left|  x\right|  ^{2n}%
\end{array}
\right)  .\label{deftildepn}%
\end{equation}
This formula shows that the degree of $A_{n}\left(  x\right)  $ is lower or
equal than $n^{2}+n.$

Let us discuss the polynomial part $\widetilde{Q}_{n}\left(  \zeta
,\theta\right)  .$ Let us write $\widetilde{P}_{n}\left(  \zeta,\theta\right)
=p_{0}\left(  \theta\right)  +p_{1}\left(  \theta\right)  \zeta+...+p_{n}%
\left(  \theta\right)  \zeta^{n}$. By formula (\ref{defpn}) it is clear that
$p_{j}\left(  \theta\right)  $ can be defined by the determinant of the matrix
in (\ref{defpn}) where we have deleted the $j$-column and the last row. An
analysis analog to the above shows that there exists a polynomial
$R_{j}\left(  x\right)  $ such that $R_{j}\left(  \zeta\theta\right)
=\zeta^{n^{2}-j}p_{j}\left(  \theta\right)  .$ Now formula (\ref{ppart}) shows
that
\[
\zeta^{n^{2}}\widetilde{Q}_{n}\left(  \zeta,\theta\right)  =\zeta^{n^{2}}%
\sum_{k=0}^{n-1}\zeta^{k}\sum_{l=0}^{n-1-k}f_{l}\left(  \theta\right)
p_{k+1+l}\left(  \theta\right)  .
\]
Since $\zeta^{n^{2}}\zeta^{k}f_{l}\left(  \theta\right)  p_{k+1+l}\left(
\theta\right)  =\zeta^{2k+1}R_{k+1+l}\left(  \zeta\theta\right)  F_{l}\left(
\zeta\theta\right)  $ one can conclude that $\frac{1}{\zeta}\zeta^{n^{2}%
}\widetilde{Q}_{n}\left(  \zeta,\theta\right)  $ is a polynomial. \endproof

We want to relate the Pad\'{e} approximation in Theorem \ref{ThmPad} to
Pad\'{e} approximation in the context of polynomials in several real
variables. Let $F_{l}$ be the homogeneous polynomial of degree $l$ defined in
Proposition \ref{P8}. The asymptotic expansion
\[
\widehat{\mu}\left(  \zeta,\theta\right)  =\sum_{l=0}^{\infty}f_{l}\left(
\theta\right)  \frac{1}{\zeta^{l+1}}=\sum_{l=0}^{\infty}F_{l}\left(
\zeta\theta\right)  \frac{1}{\zeta^{2l+1}}%
\]
and the identity $\widehat{\mu}_{\text{real}}\left(  \zeta\theta\right)
=\zeta\widehat{\mu}\left(  \zeta,\theta\right)  $, see (\ref{eqmureell}),
yield the asymptotic expansion of the real Markov transform $\widehat{\mu
}_{\text{real}}\left(  y\right)  ,$ namely
\[
\widehat{\mu}_{\text{real}}\left(  y\right)  =\sum_{l=0}^{\infty}F_{l}\left(
y\right)  \frac{1}{\left|  y\right|  ^{2l}}.
\]
By formula (\ref{eqPade}), Theorem \ref{ThmPad} and \ref{ThmPoly} we can find
polynomials $A_{n}\left(  y\right)  $ and $B_{n}\left(  y\right)  $ such that
\[
A_{n}\left(  \zeta\theta\right)  \widehat{\mu}\left(  \zeta,\theta\right)
-\zeta B_{n}\left(  \zeta\theta\right)  =\zeta^{n^{2}}\sum_{l=n}^{\infty}%
f_{l}\left(  \theta\right)  \frac{1}{\zeta^{l+1}}%
\]
for all $\theta\in\mathbb{S}^{d-1}$ such that the index $n$ is normal, i.e.
$H_{n}\left(  \mu,\theta\right)  \neq0.$ We multiply this equation by $\zeta$
and write
\[
A_{n}\left(  \zeta\theta\right)  \zeta\widehat{\mu}\left(  \zeta
,\theta\right)  -\zeta^{2}B_{n}\left(  \zeta\theta\right)  =\zeta^{n^{2}}%
\sum_{l=n}^{\infty}f_{l}\left(  \zeta\theta\right)  \frac{1}{\zeta^{2l}}.
\]
Further for the polynomial $h\left(  y\right)  =\left|  y\right|  ^{2}$ we
have $h\left(  \zeta\theta\right)  =\zeta^{2},$ so the last equation implies
for real $y=\rho\theta$ with $\left|  y\right|  >R$ and $H_{n}\left(
\mu,\theta\right)  \neq0$
\[
A_{n}\left(  y\right)  \widehat{\mu}_{\text{real}}\left(  y\right)  -\left|
y\right|  ^{2}B_{n}\left(  y\right)  =\left|  y\right|  ^{n^{2}}\sum
_{l=n}^{\infty}F_{l}\left(  y\right)  \frac{1}{\left|  y\right|  ^{2l}}.
\]
Here $A_{n}\left(  y\right)  $ and $B_{n}\left(  y\right)  $ are subject to
the conditions expressed in (\ref{eqPolAB}), and it seems to be rather
technical to convert this in direct conditions for $A_{n},B_{n}.$

We refer to \cite{Cuyt83} and \cite{CDL96} for multivariate Pad\'{e}
approximation based on polynomials in several variables.

Now we want to give an example of a measure $\mu$ such that the polynomial
$\zeta\longmapsto\widetilde{P}_{n}\left(  \zeta,\theta\right)  $ (defined in
(\ref{defpn})) is the zero polynomial. We recall at first the following result
from \cite{KoReHirosh}.

\begin{proposition}
\label{Prop6}Let $\sigma$ be a measure on $\mathbb{R}$ with compact support,
$\delta_{0}$ be the Dirac measure on $\mathbb{R}$ at the point $0$ and let
$\mu=\sigma\otimes\delta_{0}$ be the product measure. Then the multivariate
Markov transform $\widehat{\mu}$ is given by
\begin{equation}
\widehat{\sigma\otimes\delta_{0}}\left(  \zeta,e^{it}\right)  =\frac{1}%
{\omega_{2}}\sum_{l=0}^{\infty}\int x^{l}d\sigma\left(  x\right)  \frac
{\sin\left(  l+1\right)  t}{\sin t}\frac{1}{\zeta^{l+1}}\label{sigdelt}%
\end{equation}
where $\omega_{2}$ is the area measure of $\mathbb{S}^{2}.$
\end{proposition}

The last proposition has been used to show that there exists a measure $\mu$
with a support contained in an algebraic set such that $\zeta\widehat{\mu
}\left(  \zeta,\theta\right)  $ is not a rational function.

\begin{example}
\label{Ex1}Let $\sigma$ be the Lebesgue measure on $\left[  0,1\right]  ,$ so
$\int_{0}^{1}x^{l}d\sigma\left(  x\right)  =1/(l+1)$ and let $\mu
=\sigma\otimes\delta_{0}$ as in Proposition \ref{Prop6}. Then $f_{l}\left(
1\right)  =1$ for all $l\in\mathbb{N}_{0}$, and this implies that
$\widetilde{P}_{n}\left(  \zeta,1\right)  =0$ for all $n\geq2.$
\end{example}

\section{Rationality of the multivariate Mar\-kov trans\-form}

Recall that a function $f:\mathbb{R}^{d}\rightarrow\mathbb{C}$ is
\emph{rational} if there exist polynomials $p\left(  x\right)  $ and $q\left(
x\right)  \neq0$ with $f\left(  x\right)  =p\left(  x\right)  /q\left(
x\right)  $ for all $x$ with $q\left(  x\right)  \neq0.$

A theorem of Kronecker (Theorem 3.1 in \cite{NiSo91}) says that a necessary
and sufficient condition for a series $f\left(  \zeta\right)  $ of a single
variable $\zeta$ to be the Laurent expansion of a rational function is that
the Hankel determinants $H_{m}\left(  f\right)  $ are zero for all
sufficiently large $m.$

We have now the following analogue for the multivariate Markov transform:

\begin{theorem}
Let $\mu$ be a measure with support in $\overline{B_{R}}$. Then the following
statements are equivalent:

a) For each $\theta\in\mathbb{S}^{d-1}$ the function $\zeta\mapsto\widehat
{\mu}\left(  \zeta,\theta\right)  $ is rational.

b) There exists $n\in\mathbb{N}$ such that $H_{m}\left(  \mu,\theta\right)
\equiv0$ for all $m\geq n$ and $\theta\in\mathbb{S}^{d-1}.$

c) There exists $n\in\mathbb{N}$ such that for each $\theta\in\mathbb{S}%
^{d-1}$ the function $\zeta\mapsto\widehat{\mu}\left(  \zeta,\theta\right)  $
is rational of degree $\leq n.$

d) There exist polynomials $P\left(  x\right)  $ and $Q\left(  x\right)  $
such that for all $\theta\in\mathbb{S}^{d-1}$ and for all $\left|
\zeta\right|  >R$
\begin{equation}
P\left(  \zeta\theta\right)  \neq0\text{ and }\widehat{\mu}\left(
\zeta,\theta\right)  =\zeta\frac{Q\left(  \zeta\theta\right)  }{P\left(
\zeta\theta\right)  }\text{. }\label{hhaupteq}%
\end{equation}
\end{theorem}%

\proof
Assume $a)$ and let $d\left(  \theta\right)  $ be the degree of the rational
function $\zeta\mapsto\widehat{\mu}\left(  \zeta,\theta\right)  $ for
$\theta\in\mathbb{S}^{d-1}$ (recall that the degree of a rational function
$f=p/q$ with relatively prime polynomials is defined as $\max\left\{  \deg
p,\deg q\right\}  $, see \cite[p. 38]{NiSo91}). By Kronecker's theorem (cf.
\cite[p. 46]{NiSo91}) it follows that the Pad\'{e} approximant $\pi_{n}\left(
\zeta,\theta\right)  $ is equal to $\widehat{\mu}\left(  \zeta,\theta\right)
$ for all indices $n>d\left(  \theta\right)  $ and $H_{n}\left(  \mu
,\theta\right)  =0$ for all $n>d\left(  \theta\right)  .$ It follows that
$\mathbb{S}^{d-1}$ is the union of the following sets
\[
A_{n}:=\cap_{m=n}^{\infty}\left\{  \theta\in\mathbb{S}^{d-1}:H_{m}\left(
\mu,\theta\right)  =0\right\}
\]
for $n\in\mathbb{N}.$ Moreover $A_{n}$ is closed by continuity of
$\theta\mapsto H_{n}\left(  \mu,\theta\right)  $. By Baire's category theorem
there exists an index $n$ such that $A_{n}$ contains an interior point. Hence
there exists $\theta_{0}\in\mathbb{S}^{d-1}$ and a neighborhood $U$ of
$\theta_{0}$ such that $H_{m}\left(  \mu,\theta\right)  =0$ for all $\theta\in
U$ and for all $m\geq n.$ Since $\zeta^{n\left(  n-1\right)  }H_{m}\left(
\mu,\theta\right)  =\widetilde{H}_{n}\left(  \zeta\theta\right)  $ by
Proposition \ref{PropHD} we see that the polynomial $\widetilde{H}_{n}$
vanishes in a neighborhood of $\theta_{0}\in\mathbb{R}^{d}.$ Thus
$\widetilde{H}_{n}\left(  x\right)  $ is the zero polynomial and $H_{m}\left(
\mu,\theta\right)  =0$ for $\theta\in\mathbb{S}^{d-1}$ and for all $m\geq n.$
Hence we have proved $b).$

The implication $b)\rightarrow c)$ follows from Kronecker's theorem. The
implication is $c)\rightarrow a)$ is trivial.

Clearly $d)$ implies $a),$ and it suffices to show that $b)$ implies $d).$ Now
suppose that $H_{m}\left(  \mu,\theta\right)  \equiv0$ for all $m>n$ and that
$H_{n}\left(  \mu,\theta\right)  $ is not the zero function. Let $N_{n}$ be
the set of all $\theta\in\mathbb{S}^{d-1}$ such that $H_{n}\left(  \mu
,\theta\right)  \neq0,$ so $n$ is a normal index for each $\theta\in N_{n}$.
Clearly $N_{n}$ is an open non-empty set. Let $m>n$ be arbitrary. By
Proposition 3.3 in \cite[p. 45]{NiSo91} $\pi_{m}\left(  \zeta,\theta\right)
=\pi_{n}\left(  \zeta,\theta\right)  $ for all $\theta\in N_{n}$ and this
implies $\pi_{n}\left(  \zeta,\theta\right)  =\widehat{\mu}\left(
\zeta,\theta\right)  $ for all $\left|  \zeta\right|  >R$ and $\theta\in
N_{n}$. Let $\widetilde{P}_{n}\left(  \zeta,\theta\right)  $ and
$\widetilde{Q}_{n}\left(  \zeta,\theta\right)  $ as in the last section. Since
the index $n$ is normal for $\theta\in N_{n}$ it is clear that $(\widetilde
{P}_{n}\left(  \zeta,\theta\right)  ,\widetilde{Q}_{n}\left(  \zeta
,\theta\right)  )$ is an $n$-th Pad\'{e} pair and we infer
\begin{equation}
\widetilde{P}_{n}\left(  \zeta,\theta\right)  \widehat{\mu}\left(
\zeta,\theta\right)  =\widetilde{Q}_{n}\left(  \zeta,\theta\right)
.\label{schoengl}%
\end{equation}
By Theorem \ref{ThmPoly} there exist polynomials $A\left(  x\right)  ,B\left(
x\right)  $ such that
\[
A\left(  \zeta\theta\right)  =\zeta^{n^{2}}\widetilde{P}_{n}\left(
\zeta,\theta\right)  \text{ and }\zeta B\left(  \zeta\theta\right)
=\zeta^{n^{2}}\widetilde{Q}_{n}\left(  \zeta,\theta\right)
\]
for all $\zeta$ and $\theta\in\mathbb{S}^{d-1}.$ So we have for all
$y=\rho\theta$ with $\rho>R$ and $\theta\in N_{n}$ that
\begin{equation}
A\left(  y\right)  \cdot\widehat{\mu}_{\text{real}}\left(  y\right)  =\left|
y\right|  ^{2}B\left(  y\right) \label{eqid}%
\end{equation}
where $\widehat{\mu}_{\text{real}}$ is defined in (\ref{muschlange}). Since
$\widehat{\mu}_{\text{real}}\left(  y\right)  $, and obviously $A\left(
y\right)  $ and $B\left(  y\right)  ,$ are real-analytic for all
$y\in\mathbb{R}^{d},\left|  y\right|  >R$ , the equality (\ref{eqid}), valid
on an open subset of $\mathbb{R}^{d}\setminus\overline{B_{R}}$, holds for all
$y\in\mathbb{R}^{d},\left|  y\right|  >R$ as well. Moreover $\widehat{\mu
}_{\text{real}}$ has an holomorphic continuation for all $z\in\mathbb{C}^{d}$
with $L_{-}\left(  z\right)  >R$ where $L_{-}$ is defined in (\ref{defLminus}%
). So we obtain
\begin{equation}
A\left(  z\right)  \widehat{\mu}_{\text{real}}\left(  z\right)  =\left(
z_{1}^{2}+...+z_{d}^{2}\right)  B\left(  z\right)  \text{ for all }%
z\in\mathbb{C}^{d},L_{-}\left(  z\right)  >R.\label{eqAB}%
\end{equation}
Suppose now that $A$ contains an irreducible factor $g$ which has a zero
$z_{0}\in\mathbb{C}^{d}$ with $L_{-}\left(  z_{0}\right)  >R.$ By continuity
of $L_{-}$ there exists a neighborhood $U$ of $z_{0}$ with $L_{-}\left(
z\right)  >R$ for all $z\in U$. Equation (\ref{eqAB}) shows that $U\cap
g^{-1}\left\{  0\right\}  \subset B^{-1}\left\{  0\right\}  $ (recalling that
$z_{1}^{2}+...+z_{d}^{2}\neq0$ for all $z\in\mathbb{C}^{d}$ with $L_{-}\left(
z\right)  >R).$ It follows that $g$ must divide $B,$ see \cite[p. 26]{Walk}.
Inductively, we can factor out each irreducible factor of $A$ which has zero
$z_{0}\in\mathbb{C}^{d}$ with $L_{-}\left(  z_{0}\right)  >R.$ Finally we
obtain polynomials $A_{1}\left(  z\right)  $ and $B_{1}\left(  z\right)  $
such that
\[
A_{1}\left(  z\right)  \widehat{\mu}_{\text{real}}\left(  z\right)  =\left(
z_{1}^{2}+...+z_{d}^{2}\right)  B_{1}\left(  z\right)  \text{ for all }%
z\in\mathbb{C}^{d}\text{ with }L_{-}\left(  z\right)  >R,
\]
and $A_{1}\left(  z\right)  \neq0$ for all $L_{-}\left(  z_{0}\right)  >R.$
The proof is accomplished.%
\endproof

\section{A cubature formula}

Let $\sigma$ be a measure with finite moments on $\mathbb{R}$ and support in
$\left[  -R,R\right]  $ and consider the functional
\begin{equation}
T\left(  u\right)  :=\frac{1}{2\pi i}\int_{\Gamma_{R_{1}}}u\left(
\zeta\right)  \widehat{\sigma}\left(  \zeta\right)  d\zeta\label{neuComp}%
\end{equation}
where $\Gamma_{R_{1}}\left(  t\right)  =R_{1}e^{it}$ for $t\in\left[
0,2\pi\right]  $ for any $R_{1}>R$. For a polynomial $u\left(  \zeta\right)
=u_{0}+u_{1}\zeta+...+u_{m}\zeta^{m}$ we have
\begin{equation}
T\left(  u\right)  =\sum_{l=0}^{m}u_{l}\cdot f_{l}=\int_{-R}^{R}u\left(
x\right)  d\sigma\left(  x\right)  ,\label{neuComp2}%
\end{equation}
where $f_{l}:=\int_{a}^{b}x^{l}d\sigma\left(  x\right)  $ are the coefficients
in the asymptotic expansion of $\widehat{\sigma}$ given in (\ref{asympst}).

We shall make use of the following classical fact (see e.g. \cite{NiSo91}):
Let $(Q_{n},P_{n})$ be the $n$-th Pad\'{e} pair of the Markov transform
$\widehat{\sigma}\left(  \zeta\right)  $ of a non-negative measure $\sigma$
with support in the interval $[a,b] . $ If $n$ is normal (so the Hankel
determinant $H_{n}\left(  \widehat{\sigma}\right)  $ is not zero) then the
zeros $x_{1},...,x_{n}$ of $P_{n}$ are real and simple and lie in the interval
$( a,b ) ;$ moreover there exist positive coefficients $\alpha_{1}%
,...,\alpha_{n}$ such that the discrete measure
\[
\sigma_{n}=\alpha_{1}\delta_{x_{1}}+...+\alpha_{n}\delta_{x_{n}}%
\]
is identical to $\sigma$ on the subspace of all polynomials $p\left(
x\right)  $ of degree $\leq2n-1$ and we have $\alpha_{k}=Q_{n}\left(
x_{k}\right)  /P_{n}^{\prime}\left(  x_{k}\right)  $ for $k=1,...,n.$ For any
polynomial $u\left(  x\right)  $ Cauchy's theorem yields
\begin{equation}
\frac{1}{2\pi i}\int_{\Gamma_{R_{1}}}u\left(  \zeta\right)  \frac{Q_{n}\left(
\zeta\right)  }{P_{n}\left(  \zeta\right)  }d\zeta=\sum_{k=1}^{n}\alpha
_{k}u\left(  x_{k}\right)  =\int_{a}^{b}u\left(  x\right)  d\sigma
_{n}.\label{eqquad}%
\end{equation}
Combining this with (\ref{neuComp2}) we obtain the following formula
\begin{equation}
\frac{1}{2\pi i}\int_{\Gamma_{R_{1}}}u\left(  \zeta\right)  \frac{Q_{n}\left(
\zeta\right)  }{P_{n}\left(  \zeta\right)  }d\zeta=\sum_{l=0}^{2n-1}u_{l}\cdot
f_{l}.\label{eqmainpad}%
\end{equation}
valid for any polynomial $u$ of degree $\leq$ $2n-1.$

Let $\mathcal{P}\left(  \mathbb{R}^{d}\right)  $ be the set of all polynomials
in $d$ variables.

\begin{definition}
A functional $T:\mathcal{P}\left(  \mathbb{R}^{d}\right)  \rightarrow
\mathbb{C}$ is positive definite if
\[
T\left(  u^{\ast}u\right)  \geq0
\]
for all $u\in\mathcal{P}\left(  \mathbb{R}^{d}\right)  $ where $u^{\ast}$ is
the polynomial obtained from $u$ by conjugating the coefficients.
\end{definition}

\begin{definition}
A measure $\mu$ on $\mathbb{R}^{d}$ with support in the closed ball
$\overline{B_{R}}$ is called Hankel-positive if the Hankel determinants are
strictly positive, i.e.
\[
H_{n}\left(  \mu,\theta\right)  >0\text{ for all }n\in\mathbb{N},\theta
\in\mathbb{S}^{d-1}.
\]
\end{definition}

Obviously, an equivalent formulation for Hankel positivity is the requirement
that
\[
\left(  f_{l}\left(  \theta\right)  \right)  _{l=0,1,..}\text{ is strictly
positive definite}%
\]
for each $\theta\in\mathbb{S}^{d-1}.$ This means that for each $\theta
\in\mathbb{S}^{d-1}$ and for all $(x_{0},...,x_{n})\in\mathbb{R}^{n+1}%
,(x_{0},...,x_{n})\neq0$
\[
\sum_{i=0}^{n}\sum_{j=0}^{n}f_{i+j}\left(  \theta\right)  x_{i}x_{j}>0.
\]

The following result is needed in the next theorem:

\begin{proposition}
\label{ThmNozeros}Let $\widetilde{P}_{n}\left(  \zeta,\theta\right)  $ be
defined in (\ref{defpn}). If the Hankel determinant $\theta\mapsto
H_{n}\left(  \mu,\theta\right)  $ has no zeros then there exists $R_{1}>0$
such that
\begin{equation}
\widetilde{P}_{n}\left(  \zeta,\theta\right)  \neq0\text{ for all }\theta
\in\mathbb{S}^{d-1},\zeta\in\mathbb{C}\text{ with }\left|  \zeta\right|  \geq
R_{1}.\label{eqnonzero2}%
\end{equation}
\end{proposition}%

\proof
By assumption $H_{n}\left(  \mu,\theta\right)  \neq0$ for all $\theta
\in\mathbb{S}^{d-1},$ so it follows that $\zeta\longmapsto\widetilde{P}%
_{n}\left(  \zeta,\theta\right)  $ defined in (\ref{defpn}) is a polynomial of
degree exactly $n.$ Let us write
\[
\widetilde{P}_{n}\left(  \zeta,\theta\right)  =p_{0}\left(  \theta\right)
+p_{1}\left(  \theta\right)  \zeta+...+p_{n}\left(  \theta\right)  \zeta^{n}.
\]
Then $p_{n}\left(  \theta\right)  \neq0$ for all $\theta\in\mathbb{S}^{d-1} $
and $p_{n}$ is continuous. A straighforward estimate now shows that there
exists $R_{1}>0$ such that $\widetilde{P}_{n}\left(  \zeta,\theta\right)
\neq0$ for all $\left|  \zeta\right|  >R_{1}$ and for all $\theta\in
\mathbb{S}^{d-1}.$
\endproof

The following is an analog of (\ref{neuComp2}) for the multivariate Markov
transform, for the proof we refer to \cite{KoReHirosh}.

\begin{proposition}
\label{P1}\label{stielt2}Let $\mu$ be a measure on $\mathbb{R}^{d}$ with
support in $\overline{B_{R}}$ and let $R_{1}>R$. Then for any $u\in
\mathcal{P}\left(  \mathbb{R}^{d}\right)  $
\begin{equation}
M\left(  u\right)  :=\frac{1}{2\pi i\omega_{d}}\int_{\Gamma_{R_{1}}}%
\int_{\mathbb{S}^{d-1}}u\left(  \zeta\theta\right)  \widehat{\mu}\left(
\zeta,\theta\right)  d\zeta d\theta=\int_{\mathbb{R}^{d}}u\left(  x\right)
d\mu\left(  x\right)  .\label{MMintegral}%
\end{equation}
\end{proposition}

The following result is an extension of the Gau\ss\ quadrature formula to the
multivariate setting. It can be seen as a solution of the \emph{truncated
moment problem} for the class of Hankel-positive measures. We refer to
\cite{BCR84} and \cite{Fugl83} for a description of the multivariate moment problem.

\begin{theorem}
\label{T10}Let $\mu$ be a Hankel-positive measure with support in
$\overline{B_{R}}$. Let $\widetilde{P}_{n}\left(  \zeta,\theta\right)  $ and
$\widetilde{Q}_{n}\left(  \zeta,\theta\right)  $ be defined in (\ref{defpn})
and (\ref{ppart}), and let $R_{1}>R$ so large such that (\ref{eqnonzero2})
holds. Then the functional $T_{n}:\mathcal{P}\left(  \mathbb{R}^{d}\right)
\rightarrow\mathbb{C,}\ $defined by
\begin{equation}
T_{n}\left(  u\right)  :=\frac{1}{2\pi i\omega_{d}}\int_{\Gamma_{R_{1}}}%
\int_{\mathbb{S}^{d-1}}u\left(  \zeta\theta\right)  \ \frac{\widetilde{Q}%
_{n}\left(  \zeta,\theta\right)  }{\widetilde{P}_{n}\left(  \zeta
,\theta\right)  }d\zeta\ d\theta\label{defTN}%
\end{equation}
for all $u\in\mathcal{P}\left(  \mathbb{R}^{d}\right)  ,$ is positive definite
and for each polynomial $u\left(  x\right)  $ of degree $\leq2n-1$
\[
T_{n}\left(  u\right)  =\int u\left(  x\right)  d\mu\left(  x\right)  .
\]
Moreover there exists a non-negative measure $\mu_{n}$ with support in an
algebraic bounded set in $\mathbb{R}^{d}$ such that
\[
T_{n}\left(  u\right)  =\int u\left(  x\right)  d\mu_{n}\left(  x\right)
\]
for any polynomial $u.$
\end{theorem}%

\proof
1. Since (\ref{eqnonzero2}) holds for $R_{1}>R$ the expression (\ref{defTN})
is well defined. Moreover $(\widetilde{Q}_{n}\left(  \zeta,\theta\right)
,\widetilde{P}_{n}\left(  \zeta,\theta\right)  )$ is an $n$-th Pad\'{e} pair
since $H_{n}\left(  \mu,\theta\right)  \neq0$ for all $\theta\in
\mathbb{S}^{d-1}$ by Hankel positivity.

2. Suppose that $u$ is a function of the form
\begin{equation}
u\left(  \zeta\theta\right)  =u_{0}\left(  \theta\right)  +...+u_{2n-1}\left(
\theta\right)  \zeta^{2n-1}\label{neuuform}%
\end{equation}
with continuous functions $u_{0},...,u_{2n-1}.$ By (\ref{eqmainpad}) we have
\begin{equation}
\frac{1}{2\pi i}\int_{\Gamma_{R_{1}}}u\left(  \zeta\theta\right)  \frac
{Q_{n}\left(  \zeta,\theta\right)  }{P_{n}\left(  \zeta,\theta\right)  }%
d\zeta=\sum_{l=0}^{2n-1}u_{l}\left(  \theta\right)  f_{l}\left(
\theta\right)  .\label{eqpad2}%
\end{equation}
Integration over $\mathbb{S}^{d-1}$ gives
\begin{equation}
T_{n}\left(  u\right)  =\frac{1}{\omega_{d}}\sum_{l=0}^{2n-1}\int
_{\mathbb{S}^{d-1}}u_{l}\left(  \theta\right)  f_{l}\left(  \theta\right)
d\theta.\label{neuTform}%
\end{equation}

3. Let $u$ be a polynomial. Proposition \ref{P1} shows that
\begin{equation}
\int u\left(  x\right)  d\mu\left(  x\right)  =\frac{1}{2\pi i\omega_{d}}%
\int_{\mathbb{S}^{d-1}}\int_{\Gamma_{R_{1}}}u\left(  \zeta\theta\right)
\widehat{\mu}\left(  \zeta,\theta\right)  d\zeta d\theta.\label{uformel}%
\end{equation}
If $u$ has degree $\leq2n-1$ then for each $\theta\in\mathbb{S}^{d-1}$ the
function $\zeta\mapsto u\left(  \zeta\theta\right)  $ is a polynomial of
degree $\leq2n-1,$ so $u$ is of the form (\ref{neuuform}). Insert
(\ref{eqmuasym}) in (\ref{uformel}) and integrate over $\Gamma_{R_{1}}$ to
obtain
\begin{equation}
\int u\left(  x\right)  d\mu\left(  x\right)  =\frac{1}{\omega_{d}}\sum
_{l=0}^{2n-1}\int_{\mathbb{S}^{d-1}}u_{l}\left(  \theta\right)  f_{l}\left(
\theta\right)  d\theta.\label{muflgleich}%
\end{equation}
Comparing (\ref{neuTform}) with (\ref{muflgleich}) we conclude that
$T_{n}\left(  u\right)  $ is equal to $\int u\left(  x\right)  d\mu\left(
x\right)  $ for any polynomial of degree $\leq2n-1.$

4. Let us discuss the question of positive definiteness of $T_{n}$. Let
$R\left(  x\right)  $ be a real-valued polynomial. We have to show that
$T_{n}\left(  R^{2}\right)  \geq0.$ By the euclidean algorithm applied to the
polynomials $\zeta\mapsto R\left(  \zeta\theta\right)  $ and $\zeta
\mapsto\widetilde{P}_{n}\left(  \zeta\theta\right)  $ for each fixed $\theta,$
there exist a polynomial $\zeta\mapsto d\left(  \zeta,\theta\right)  ,$ and a
polynomial $\zeta\mapsto e\left(  \zeta,\theta\right)  $ of degree $<n,$ such
that
\[
R\left(  \zeta\theta\right)  =d\left(  \zeta,\theta\right)  \widetilde{P}%
_{n}\left(  \zeta,\theta\right)  +e\left(  \zeta,\theta\right)  .
\]
Write $e\left(  \zeta,\theta\right)  =e_{0}\left(  \theta\right)
+...+e_{n-1}\left(  \theta\right)  \zeta^{n-1}.$ Then
\[
\left(  R\left(  \zeta\theta\right)  \right)  ^{2}=d^{2}\left(  \zeta
,\theta\right)  \left(  \widetilde{P}_{n}\left(  \zeta,\theta\right)  \right)
^{2}+2d\left(  \zeta,\theta\right)  e\left(  \zeta,\theta\right)
\widetilde{P}_{n}\left(  \zeta,\theta\right)  +e^{2}\left(  \zeta
,\theta\right)  .
\]
Multiply the last equation with $\widetilde{Q}_{n}\left(  \zeta,\theta\right)
/\widetilde{P}_{n}\left(  \zeta,\theta\right)  $ and integrate with respect to
$\zeta$ over $\Gamma_{R_{1}}$. Then
\begin{equation}
b\left(  \theta\right)  :=\frac{1}{2\pi i}\int_{\Gamma_{R_{1}}}R^{2}\left(
\zeta\theta\right)  \frac{\widetilde{Q}_{n}\left(  \zeta,\theta\right)
}{\widetilde{P}_{n}\left(  \zeta,\theta\right)  }d\zeta=\frac{1}{2\pi i}%
\int_{\Gamma_{R_{1}}}e^{2}\left(  \zeta,\theta\right)  \frac{\widetilde{Q}%
_{n}\left(  \zeta,\theta\right)  }{\widetilde{P}_{n}\left(  \zeta
,\theta\right)  }d\zeta.\label{bteta}%
\end{equation}
Since $\zeta\rightarrow e^{2}\left(  \zeta,\theta\right)  $ is a polynomial of
degree $\leq2n-1$, (\ref{eqpad2}) yields
\begin{equation}
b\left(  \theta\right)  =\frac{1}{2\pi i}\int_{\Gamma_{R_{1}}}e^{2}\left(
\zeta,\theta\right)  \frac{\widetilde{Q}_{n}\left(  \zeta,\theta\right)
}{\widetilde{P}_{n}\left(  \zeta,\theta\right)  }d\zeta=\sum_{k,l=0}%
^{n-1}e_{k}\left(  \theta\right)  e_{l}\left(  \theta\right)  f_{k+l}\left(
\theta\right)  .\label{bteta2}%
\end{equation}
Integrate the last equation with respect to $\mathbb{S}^{d-1}$ and use the
definition of $T_{n}$ in order to obtain
\[
T_{n}\left(  R^{2}\right)  =\frac{1}{\omega_{d}}\int_{\mathbb{S}^{d-1}}%
(\sum_{k,l=0}^{n-1}e_{k}\left(  \theta\right)  e_{l}\left(  \theta\right)
f_{k+l}\left(  \theta\right)  )d\theta.
\]
Since $\left(  f_{l}\right)  _{l}$ is strictly positive definite we know that
$\sum_{k,l=0}^{n-1}e_{k}\left(  \theta\right)  e_{l}\left(  \theta\right)
f_{k+l}\left(  \theta\right)  \geq0$ for each $\theta\in\mathbb{S}^{d-1}$, in
particular $T_{n}\left(  R^{2}\right)  \geq0.$

6. Let $A_{n}$ be the polynomial defined in Theorem \ref{ThmPoly} such that
$A_{n}\left(  \zeta\theta\right)  =\zeta^{n^{2}}\widetilde{P}_{n}\left(
\zeta\theta\right)  .$ Note that $A_{n}$ has real coefficients. It follows
that for any polynomial $u$
\[
T_{n}\left(  A_{n}u\right)  =\frac{1}{2\pi i\omega_{d}}\int_{\mathbb{S}^{d-1}%
}\int_{\Gamma_{R_{1}}}\zeta^{n^{2}}u\left(  \zeta\theta\right)  \widetilde
{Q}_{n}\left(  \zeta,\theta\right)  d\zeta\ d\theta=0
\]
since $\zeta\longmapsto\zeta^{n^{2}}u\left(  \zeta\theta\right)  \widetilde
{Q}_{n}\left(  \zeta,\theta\right)  $ is a polynomial. The polynomial
$\zeta\longmapsto\widetilde{P}_{n}\left(  \zeta\theta\right)  $ has only zeros
in the interval $\left(  -R,R\right)  ,$ so it follows that $\widetilde{P}%
_{n}\left(  \rho\theta\right)  \neq0$ for all $\rho>R.$ Hence $A_{n}\left(
\rho\theta\right)  \neq0$ for all $\rho>R,$ so the zero set of $y\longmapsto
A_{n}\left(  y\right)  $ is contained in the ball $\overline{B_{R}}.$

7. The existence of a representation measure $\mu_{n}$ follows from Theorem
\ref{Schmudresult} below, cf. \cite{Schm91a} (applied to $m=2$ and
$f_{1}=A_{n}$ and $f_{2}=-A_{n}).$ The proof is finished.%
\endproof

\begin{theorem}
\label{Schmudresult}(\textbf{Schm\"{u}dgen}) Let $\frak{S}$ $:\mathcal{P}%
\left(  \mathbb{R}^{d}\right)  \rightarrow\mathbb{C}$ be a positive definite
functional and let $f_{1},...,f_{m}$ be polynomials with real coefficients
such that the set
\[
K:=\left\{  x\in\mathbb{R}^{d}:f_{j}\left(  x\right)  \geq0\text{ for all
}j=1,2,...,m\right\}
\]
is compact. Then there exists a measure $\mu$ with support in $K$ representing
$\frak{S}$ if and only if
\begin{equation}
\frak{S}\left(  f_{j_{1}}...f_{j_{s}}\cdot p^{\ast}p\right)  \geq
0\label{rrrcondition}%
\end{equation}
for all pairwise different $j_{1},...,j_{s}\in\left\{  1,...,m\right\}  $ and
for all $p\in\mathcal{P}\left(  \mathbb{R}^{d}\right)  .$
\end{theorem}

It is not difficult to see that a rotation invariant measure with
non-algebraic support is Hankel-positive, cf. Theorem \ref{Rotthm}. Now we
give a different class of examples:

\begin{proposition}
Let $w_{0},w_{1}:\left[  0,\infty\right)  \rightarrow\mathbb{R}$ be bounded
continuous functions with compact support such that $\left|  w_{1}\left(
r\right)  \right|  \leq w_{0}\left(  r\right)  $ for all $r\geq0,$ and
$w_{0}\neq0.$ Assume that the measure $\mu$ has the density $d\mu:=w\left(
r,\vartheta\right)  rdrd\vartheta$ where
\[
w\left(  r,\vartheta\right)  =w_{0}\left(  r\right)  +w_{1}\left(  r\right)
\cos\vartheta\text{ }%
\]
for all $r>0,\vartheta\in\left[  0,2\pi\right]  .$ Then $\mu$ is Hankel-positive.
\end{proposition}%

\proof
Note that the assumption $\left|  w_{1}\left(  r\right)  \right|  \leq
w_{0}\left(  r\right)  $ for all $r\geq0$ assures that
\begin{equation}
w\left(  r,\vartheta\right)  =w_{0}\left(  r\right)  +w_{1}\left(  r\right)
\cos\vartheta\geq0\label{eqcondF}%
\end{equation}
for all $r>0,\vartheta\in\left[  0,2\pi\right]  .$ Let us write $\theta
=e^{i\vartheta}$ with $\vartheta\in\left[  0,2\pi\right)  .$ Note that
$Y_{0}\left(  \theta\right)  =1/\sqrt{2\pi}$ and $Y_{k,1}\left(
\theta\right)  =\frac{1}{\sqrt{\pi}}\cos k\vartheta$ and $Y_{k,2}\left(
\theta\right)  =\frac{1}{\sqrt{\pi}}\sin k\vartheta,$ $k\in\mathbb{N}_{0}$,
provides an orthonormal basis of spherical harmonics. Then
\[
\int_{\mathbb{R}^{d}}\left|  x\right|  ^{2s}Y_{0}\left(  x\right)  d\mu
=\frac{1}{\sqrt{2\pi}}\int_{0}^{\infty}\int_{0}^{2\pi}r^{2s}w\left(
r,\vartheta\right)  rdrd\vartheta=\sqrt{2\pi}\int_{0}^{\infty}r^{2s+1}%
w_{0}\left(  r\right)  dr,
\]
and
\begin{align*}
\int_{\mathbb{R}^{d}}\left|  x\right|  ^{2s}Y_{1,1}\left(  x\right)  d\mu &
=\frac{1}{\sqrt{\pi}}\int_{0}^{\infty}\int_{0}^{2\pi}r^{2s}\cdot
r\cos\vartheta\cdot w\left(  r,\vartheta\right)  rdrd\vartheta\\
& =\sqrt{\pi}\int_{0}^{\infty}r^{2s+2}w_{1}\left(  r\right)  dr,
\end{align*}
while all other distributed moments $c_{s,k,m}$ are zero. By Theorem \ref{T3}
we obtain the Markov transform:
\[
\widehat{\mu}\left(  \zeta,e^{i\vartheta}\right)  =\sum_{s=0}^{\infty}\frac
{1}{\zeta^{2s+1}}\int_{0}^{\infty}r^{2s+1}w_{0}\left(  r\right)  dr+\sum
_{s=0}^{\infty}\frac{\cos\vartheta}{\zeta^{2s+2}}\int_{0}^{\infty}%
r^{2s+2}w_{1}\left(  r\right)  dr.
\]
So $f_{2s}\left(  e^{i\vartheta}\right)  =\int_{0}^{\infty}r^{2s+1}%
w_{0}\left(  r\right)  dr$ and $f_{2s+1}\left(  e^{i\vartheta}\right)
=\cos\vartheta\int_{0}^{\infty}r^{2s+2}w_{1}\left(  r\right)  dr.$

Extend the function $w_{0}$ to an odd function $w_{0}^{\text{odd}}$ on
$\mathbb{R}\setminus\left\{  0\right\}  $, so define $w_{0}^{\text{odd}%
}\left(  -r\right)  :=-w_{0}\left(  r\right)  $ for $r>0,$ and extend $w_{1}$
to an even function $w_{1}^{\text{ev}}$, so $w_{1}^{\text{ev}}\left(
-r\right)  =w_{1}\left(  r\right)  $ for $r>0.$ Define a function
$G_{\vartheta}:\mathbb{R}\rightarrow\mathbb{R}$ by
\[
G_{\vartheta}\left(  r\right)  :=r\cdot\left[  w_{0}^{\text{odd}}\left(
r\right)  +w_{1}^{\text{ev}}\left(  r\right)  \cos\vartheta\right]  .
\]
Note that $G_{\vartheta}\left(  r\right)  \geq0$ for all $r\geq0$ by condition
(\ref{eqcondF}). Moreover $G_{\vartheta}\left(  -r\right)  =rw_{0}\left(
r\right)  -rw_{1}\left(  r\right)  \cos\vartheta\geq0$ for $r>0 $ again by
(\ref{eqcondF}). Hence $G_{\vartheta}\left(  r\right)  \geq0$ for all
$r\in\mathbb{R}.$ A straightforward calcuation shows that
\[
f_{l}\left(  e^{i\vartheta}\right)  =\frac{1}{2}\int_{-\infty}^{\infty}%
r^{l}\ G_{\vartheta}\left(  r\right)  rdr
\]
for all $l\in\mathbb{N}_{0}.$ This shows that $\left(  f_{l}\left(
e^{i\vartheta}\right)  \right)  _{l\in\mathbb{N}_{0}}$ is a positive definite
sequence. If the sequence is not strictly positive definite then there exists
a polynomial $p\left(  r\right)  \neq0$ such that
\[
\frac{1}{2}\int_{-\infty}^{\infty}\left(  p\left(  r\right)  \right)
^{2}\ G_{\vartheta}\left(  r\right)  rdr=0.
\]
Since $G_{\vartheta}\left(  r\right)  $ is continuous on $\mathbb{R}%
\setminus\left\{  0\right\}  $ this implies that $G_{\vartheta}\left(
r\right)  =0$ for all $r\neq0.$ Then $0=G_{\vartheta}\left(  r\right)
+G_{\vartheta}\left(  -r\right)  =2rw_{0}\left(  r\right)  $ for all $r>0,$ a
contradiction to our assumption $w_{0}\neq0$
\endproof

In Theorem \ref{Rotthm} we have seen that the Markov transform $\widehat{\mu
}\left(  \zeta,\theta\right)  $ of a rotation invariant measure has the
property that $\zeta\longmapsto\widehat{\mu}\left(  \zeta,\theta\right)  $
posseses an analytic continuation to the upper half plane. Next we show that
the same is true for Hankel-positive measures.

\begin{theorem}
Let $\mu$ be a finite measure on $\mathbb{R}^{d}$ with support in
$\overline{B_{R}}$. If $\mu$ is Hankel-positive then for each $\theta
\in\mathbb{S}^{d-1}$ the function $\zeta\longmapsto\widehat{\mu}\left(
\zeta,\theta\right)  $ posseses an analytic continuation to the upper half
plane such that
\[
\text{Im}\widehat{\mu}\left(  \zeta,\theta\right)  \leq0\text{ for all
}\text{Im}\zeta>0,\theta\in\mathbb{S}^{d-1}.
\]
\end{theorem}%

\proof
Suppose that the sequence $\left(  f_{l}\left(  \theta\right)  \right)
_{l=0,1,..}$ is strictly positive definite. By the solution of the Hamburger
moment problem (p. 65 in \cite{NiSo91}) there exists a finite non-negative
measure $\sigma_{\theta}$ on $\mathbb{R}$ such that $\widehat{\mu}\left(
\zeta,\theta\right)  =\int\frac{1}{\zeta-t}d\sigma_{\theta}\left(  t\right)
.$ Hence $\zeta\longmapsto\widehat{\mu}\left(  \zeta,\theta\right)  $ extends
to the upper half plane for $\zeta$ and the condition $\text{Im}$
$\widehat{\mu}\left(  \zeta,\theta\right)  \leq0$ for all $\text{Im}\zeta>0$
and $\theta\in\mathbb{S}^{d-1}$ follows from this integral representation.%
\endproof

Note that the measure $\sigma_{\theta}$ in the last proof has the property
that its support set is infinite since $\left(  f_{l}\left(  \theta\right)
\right)  _{l=0,1,..}$ is \emph{strictly} positive definite. If we know that
$\text{Im}\widehat{\mu}\left(  \zeta,\theta\right)  \leq0$ for all $\theta
\in\mathbb{S}^{d-1}$ and for all $\text{Im}\zeta>0$ then the function
$g_{\theta}$ defined by $g_{\theta}\left(  \zeta\right)  :=\widehat{\mu
}\left(  \zeta,\theta\right)  $ is in the Nevanlinna class (see \cite{Akhi65})
and the coefficients of the Laurent expansion are exactly the numbers
$f_{l}\left(  \theta\right)  .$ Hence we can conclude that the sequence
$\left(  f_{l}\left(  \theta\right)  \right)  _{l}$ is positive definite for
each $\theta\in\mathbb{S}^{d-1}$. Note that Hankel positivity means that the
sequence $\left(  f_{l}\left(  \theta\right)  \right)  _{l}$ is strictly
positive definite for each $\theta\in\mathbb{S}^{d-1}$.

Let us remark that in \cite{KoRe06b} we have introduced a different method for
approximating a large class of signed measures, the so-called
\emph{pseudo-positive measures}, and we have provided a multivariate
generalization of the Gauss-Jacobi quadrature formula. 

We thank the anonimous referee for the careful reading of the manuscript and
for pointing to us  a number of errors.

1. Ognyan Kounchev, Institute of Mathematics and Informatics, Bulgarian
Academy of Sciences, 8 Acad. G. Bonchev Str., 1113 Sofia, Bulgaria;

e--mail: kounchev@math.bas.bg, kounchev@math.uni--duisburg.de

2. Hermann Render, Departamento de Matem\'{a}ticas y Computati\'{o}n,
Universidad de la Rioja, Edificio Vives, Luis de Ulloa, s/n. 26004
Logro\~{n}o, Spain; e-mail: render@gmx.de; hermann.render@unirioja.es
\end{document}